\documentclass[preprint,10pt]{elsarticle} 
\usepackage[english]{babel}

\usepackage{enumerate}     
\usepackage{amsmath}
\usepackage{amssymb}
\usepackage{graphicx}    
\usepackage{epsfig} 

\def\corr{\mathrm{(corr)}}
\def\vc{\zeta^{\corr}_{\epsilon}(t)}
\def\vca{\zeta^{\corr}_{\epsilon}(t_i)}
\def\vcc{\zeta^{\corr}_{\epsilon}(t_m)}

\def\st{\sqrt{k^2-4u(t)}}
\def\a{t_i}
\def\b{t_f}
\def\c{t_m}
\def\dz{\mathrm{d}z} 
\def\ds{\mathrm{d}s}
\def\dy{\mathrm{d}y}
\journal{Systems \& Control Letters (October 26, 2011)}  
\begin{document}

\makeatletter
\def\tg{\mathop{\operator@font tg}\nolimits}
\def\cotg{\mathop{\operator@font cotg}\nolimits}
\def\arcsin{\mathop{\operator@font arcsin}\nolimits}
\def\arccos{\mathop{\operator@font arccos}\nolimits}
\def\arctg{\mathop{\operator@font arctg}\nolimits} 
\def\arccotg{\mathop{\operator@font arccotg}\nolimits}
\def\ln{\mathop{\operator@font ln}\nolimits}
\makeatother

%%%%%%%%%%%%%%%%%

\begin{abstract}
A new systematic approach to the construction of approximate solutions to a class of nonlinear singularly perturbed feedback control systems using the boundary layer functions especially with regard to the possible occurrence of the boundary layers is proposed. For example, problems with feedback control, such as the steady-states of the thermostats, where the controllers add or remove heat, depending upon the temperature registered in another place of the heated bar, can be interpreted with a second-order ordinary differential equation subject to a  nonlocal three--point boundary condition.  The $O(\epsilon)$ accurate approximation of behavior of these nonlinear systems in terms of the exponentially small boundary layer functions is given. At the end of this paper, we formulate the unsolved controllability problem for nonlinear systems.
\end{abstract}


\title{On the approximation of the boundary layers for the controllability problem of nonlinear singularly perturbed systems}
\author{R.~Vrabel}
\ead{robert.vrabel@stuba.sk}
\address{Institute of Applied Informatics, Automation and Mathematics, Faculty of Materials Science and Technology, Hajdoczyho 1,  917 01 Trnava,   Slovakia}

\newtheorem{thm}{Theorem}
\newtheorem{lem}[thm]{Lemma}
\newdefinition{rmk}{Remark}
\newdefinition{ex}{Example}
\newproof{pf}{Proof}
\newproof{pot1}{Proof of Theorem \ref{thm1}}
\newproof{pot2}{Proof of Theorem \ref{thm2}}

\pagestyle{headings}

\maketitle

\section{Motivation and introduction}

In various fields of science and engineering, systems with two-time-scale dynamics are often investigated. In state space, such systems are commonly modeled using the mathematical framework of singular perturbations, with a small parameter, say $\epsilon$, determining the degree of separation between the "slow" and "fast" channels of the system. Singularly perturbed systems (SPS) can also occur due to the presence of small "parasitic" parameters, armature inductance in a common model for most DC motors, small time constants, etc.

Singular perturbation problems arise also in  heat transfer problem with large Peclet numbers (we often assume $\epsilon$ to be small in order to diminish the effect of diffusion (\cite{MuZua}), Navier-Stokes flows with large Reynolds numbers, chemical reactor theory, aerodynamics, control of reaction-diffusion processes (\cite{CaZu}, \cite{MiZu}), quantum mechanics (\cite{AlGaGlLaN}), optimal control (\cite{NG}), for example.  

The literature on control of nonlinear SPS is extensive, at least starting with the pioneering work of P. Kokotovic {\it et al.} nearly 30 years ago (\cite{KoKR}) and continuing to the present including authors such as Z. Artstein (\cite{Arn1, Arn2}), V. Gaitsgory (\cite{ArnG, Ga, GaNg}), etc (see,  e.g. \cite{BB, BEB, Chte, MCZ, MeJi} and the references therein).

\section{Problem formulation} 

In this paper, we will consider the nonlinear singularly perturbed feedback control system without an outer disturbance of the form
{\setlength\arraycolsep{1pt}
\begin{eqnarray}
y'(t)&=&w(t)\label{def_system1}  \\ 
\epsilon w'(t)&=&-ky(t)+f\left(u(t),y(t)\right)\label{def_system2}\\
v(t)&=&g(y(t))\label{def_system3}
\end{eqnarray}}
with  the required nonlocal boundary conditions 
\begin{equation}\label{def_BC}
v(\a)=v(\c)=v(\b),\quad \a<\c<\b, 
\end{equation}
where $\epsilon>0$ is a small perturbation parameter, $[y,w]^T$ is the state vector, $v(t)$ is the
measured output, $u(t)$ is the input control, $k<0$ is a constant and  $g$ is a monotone increasing (decreasing) function on $\mathbb{R}$. The state and control variables are not constrained by any boundaries, initial time $\a$ and final time $\b$ are fixed and $y(\a),$ $y(\b)$ are free. 

Such boundary value problems can arise in the study of the steady--states of a heated bar with the thermostats, where the controllers at $t=\a$ and $t=\b$ maintain a temperature according to the temperature detected by a sensor at $t=\c.$ In this case, we consider a uniform bar of length $\b-\a$ with non-uniform temperature lying on the $t$-axis from $t=\a$ to $t=\b.$ The parameter $\epsilon$ represents the thermal diffusivity.

Different from \cite{Bog}, in this paper we will not assume that $y(\a)$ and $y(\b)$ are fixed and moreover we investigate three-point boundary value problem. There have been some papers considered the  multi-point boundary value problems in the literature (see, e.g.
\cite{Jan}, \cite{Kha}, \cite{KhaWe}, \cite{Xu}) by applying the well known coincidence degree theory and Schauder fixed point theorem or the method of lower and upper solutions.  However, there have been fewer papers considered the three--point boundary value problems for SPS without the derivative in the boundary conditions. Recently, in the paper \cite{LL}, it has been studied the nonlinear system of the form $\epsilon^2y''=f(t,y,y'),$ $0<t<1$ subject to the boundary conditions $y(0)=0,$ $y(1)=py(\tau),$ $0<\tau<1$ and $p<1,$ where the assumption $p<1$ was crucial for proving the main result. 

One of the typical behaviors of SPS is the boundary layer phenomenon: the solutions vary rapidly within very thin layer regions near the boundary.
The novelty of our approach lies in the introduction of the exponentially small boundary layer functions 
into the analysis of nonlocal boundary value problems and approximation of their solutions. The situation in the case of nonlocal boundary value problem is complicated by the fact that there is an inner point in the boundary conditions, in contrast to  the "standard"  boundary conditions as the Dirichlet problem, Neumann problem, Robin problem, periodic boundary value problem (\cite{ChHo}, \cite{CoHa}), for example. In the problem considered, there does not exist a positive solution $\tilde \zeta_\epsilon$ of differential equation $\epsilon y''-my=0,$ $m>0,$ $\epsilon>0$ (that is, $\tilde \zeta_\epsilon$ is convex) such that $\tilde \zeta_\epsilon(\c)-\tilde \zeta_\epsilon(\a)=\eta(\c)-\eta(\a)>0$ and $\tilde \zeta_\epsilon(t)\rightarrow0^+$ for $t\in (\a,\b]$ and $\epsilon\rightarrow 0^+,$ which could be used to solve this problem by the method of lower and upper solutions and consequently, to approximate the solutions. The application of convex functions is essential for composing the appropriate barrier functions for two-endpoint boundary conditions, see, e.g. \cite{ChHo}.  

The following assumptions will be made throughout the paper.
\begin{itemize}
\item[{\bf A1.}] For limiting problem (in (\ref{def_system2}) letting $\epsilon\rightarrow 0^+$) $ky=f\left(u(t),y\right)$ there exists $C^2$ function $\eta=\eta(t)$ (that is, $\eta$ is continuous up to second derivative) such that $k\eta(t)=f\left(u(t),\eta(t)\right)$ on $[\a,\b].$
\end{itemize}

Denote ${H}(\eta)=\left\{ (t,y);\quad \a\leq t\leq \b,  \vert y-\eta (t)\vert <d(t)\right\},$ where $d(t)$  is the positive continuous function on $[ \a,\b]$ such that
\begin{displaymath}
d(t) = \left\{ \begin{array}{ll}
\vert \eta (\c)-\eta(\a)\vert+\delta & \textrm{for $\a\leq t\leq \a+\frac\delta2$}\\
\delta  & \textrm{for $\a+\delta\leq t\leq \b-\delta,$}\\
\vert \eta(\b)-\eta(\c)\vert+\delta & \textrm{for $\b-\frac\delta2\leq t\leq \b$}
\end{array} \right.
\end{displaymath}
$\delta$ is a small positive constant. 

\begin{itemize}
\item[{\bf A2.}] The function $f\in C^1({H}(\eta ))$ satisfies the condition
 \begin{equation*}\label{cond}
\left\vert\frac{\partial f(u(t),y)}{\partial y}\right\vert\leq \lambda <-k \quad\mathrm{for\  every}\quad (t,y)\in {H}(\eta).
\end{equation*}
\end{itemize}

The assumption (A2) means that the linearization of SPS (\ref{def_system1}), (\ref{def_system2}) in a neighbourhood of the set $[\eta(t),0],$ $t\in[\a,\b],$ as a set of critical points, has no eigenvalues on the imaginary axis.

In this paper, we characterize the dynamics for slow variable $y$ in a neighborhood of $\eta (t)$ for sufficiently small values of the singular perturbation parameter $\epsilon$ and $t\in[\a,\b].$ Especially, we focus our attention on the appearance of boundary layers. Moreover, we give the $O(\epsilon)$ accurate  approximation of $y$ on $[\a,\b].$
 
Obviously, $y$ is a solution of boundary value problem 
\begin{equation}\label{def_DE2}
\epsilon y''(t)+ky(t)=f\left(u(t),y(t)\right)
\end{equation}
\begin{equation}\label{def_BC2}
y(\a)=y(\c)=y(\b),\quad \a<\c<\b.
\end{equation}

Recently in \cite{Vra} we have shown that the solutions of (\ref{def_DE2}), (\ref{def_BC2}), in general, start with fast transient ($\left\vert w_\epsilon(\a)\right\vert\rightarrow \infty$)  of $y_\epsilon(t)$ from $y_\epsilon(\a)$ to $\eta (t),$ which is the so--called boundary layer phenomenon, and after decay of this transient they remain close to $\eta (t)$ with an arising new fast transient of $y_\epsilon(t)$ from $\eta (t)$ to $y_\epsilon(\b)$ ($\left\vert w_\epsilon(\b)\right\vert\rightarrow \infty$).  Boundary layers are formed due to the nonuniform convergence of
the exact solution $y_\epsilon$ to the degenerate solution $\eta$ in the neighborhood of the ends  $\a$ and $\b$ of the considered interval.

\section{Behavior of SPS for \lowercase{$\epsilon\rightarrow 0^+$}}
%\iffalse
\begin{thm}[compare with \cite{Vra}, Theorem 2.1]\label{maintheorem}
Under the assumptions \\ (A1) and (A2) there exists $\epsilon_0 $ such that for every $\epsilon \in (0,\epsilon_0]$ and for every input control $u$ the SPS (\ref{def_DE2}), (\ref{def_BC2}) has in ${H}(\eta)$ an unique realization, $y_\epsilon,$ satisfying the inequality 
\begin{equation*}%\label{...}
-\vc-\hat \zeta_\epsilon(t)-C\epsilon\leq y_\epsilon(t)-(\eta (t)+\zeta_\epsilon(t))\leq\hat \zeta_\epsilon(t)+ C\epsilon
\end{equation*}
for $\eta (\c)-\eta (\a)\geq 0$ and
\begin{equation*}%\label{...}
-\hat \zeta_\epsilon(t)-C\epsilon\leq y_\epsilon(t)-(\eta (t)+\zeta_\epsilon(t))\leq\vc+\hat \zeta_\epsilon(t)+ C\epsilon
\end{equation*}
for $\eta (\c)-\eta (\a)\leq 0$  on $[ \a,\b]$ where
 {\setlength\arraycolsep{1pt}
\begin{eqnarray*}
\zeta_\epsilon(t) & = & \frac{\eta (\c)-\eta (\a)}{D}\cdot \Big( e^{\sqrt{\frac m{\epsilon}}(\b-t)}-e^{\sqrt{\frac m{\epsilon}}(t-\b)} \\
&+&e^{\sqrt{\frac m{\epsilon}}(t-\c)}-e^{\sqrt{\frac m{\epsilon}}(\c-t)}\Big),\\
\hat \zeta_\epsilon(t) & = &\frac{\vert \eta (\b)-\eta (\c)\vert}{D}\cdot\Big( e^{\sqrt{\frac m{\epsilon}}(t-\a)}-e^{\sqrt{\frac m{\epsilon}}(\a-t)} \\
&+&e^{\sqrt{\frac m{\epsilon}}(\c-t)}-e^{\sqrt{\frac m{\epsilon}}(t-\c)}\Big),\\
D & = &\left(e^{\sqrt{\frac m{\epsilon}}(\b-\a)}+e^{\sqrt{\frac m{\epsilon}}(\c-\b)}+e^{\sqrt{\frac m{\epsilon}}(\a-\c)}\right)\\ 
&-&\left(e^{\sqrt{\frac m{\epsilon}}(\a-\b)}+e^{\sqrt{\frac m{\epsilon}}(\b-\c)}+e^{\sqrt{\frac m{\epsilon}}(\c-\a)}\right),
\end{eqnarray*}}
$m=-k-\lambda,$  $C=\frac1m\max\left\{\left\vert \eta'' (t)\right\vert; t\in[ \a,\b]\right\}$  and the positive function 
{\setlength\arraycolsep{1pt}
\begin{eqnarray*}
\vc&=&\frac{\lambda\vert \eta (\c)-\eta (\a)\vert}{\sqrt{m\epsilon}}\cdot\left[-{ O}(1)\frac{\zeta_\epsilon(t)}{(\eta (\c)-\eta (\a))}\right.  \\
 &+&\left.{ O}\left(e^{\sqrt{\frac{m}{\epsilon}}(\a-\c)}\right)\frac{\hat \zeta_\epsilon(t)}{\vert \eta (\b)-\eta (\c)\vert}+t{ O}\left(e^{\sqrt{\frac{m}{\epsilon}}\chi (t)}\right)\right],
\end{eqnarray*}}
$\chi (t)<0$ for $t\in(\a,\b]$ and $\vca=\vcc.$
\end{thm}

\centerline{}

We write $s(\epsilon)={O}(r(\epsilon))$ when $0<\lim\limits_{\epsilon\rightarrow0^+ }\left\vert\frac{s(\epsilon)}{ r(\epsilon)}\right\vert<\infty.$ 

\centerline{}

The function $\zeta_{\epsilon}(t)$ satisfies
\begin{enumerate}
\item $\epsilon \zeta''_{\epsilon}-m\zeta_{\epsilon}=0,$
\item $\zeta_{\epsilon}(\c)-\zeta_{\epsilon}(\a)=-(\eta (\c)-\eta (\a)),$ $\zeta_{\epsilon}(\b)-\zeta_{\epsilon}(\c)=0,$
\item $\zeta_{\epsilon}(t)\geq 0$ $(\leq 0)$ is decreasing (increasing) for $\a\leq t\leq\frac{\b+\c}{2}$ and increasing (decreasing) for $\frac{\b+\c}{2}\leq t\leq \b$  if $\eta (\c)-\eta (\a)\geq 0$ $(\leq 0),$
\item $\zeta_{\epsilon}(t)$ converges uniformly to $0$ for $\epsilon\rightarrow 0^+$ on every compact subset of $(\a, \b],$
\item $\zeta_{\epsilon}(t)=(\eta (\c)-\eta (\a)){O}\left(e^{\sqrt{\frac{m}{\epsilon}}\chi(t)}\right)$ where $\chi (t)=\a-t$ for $\a\leq t\leq\frac{\b+\c}{2}$ and  $\chi (t)=t-\b+\a-\c$ for $\frac{\b+\c}{2}<t\leq \b.$
\end{enumerate}

\centerline{}

The function $\hat \zeta_{\epsilon}(t)$ satisfies
\begin{enumerate}
\item $\epsilon\hat \zeta_\epsilon''-m\hat \zeta_\epsilon=0,$
\item $\hat \zeta_{\epsilon}(\c)-\hat \zeta_{\epsilon}(\a)=0,$ $\hat \zeta_{\epsilon}(\b)-\hat \zeta_{\epsilon}(\c)=\vert \eta (\b)-\eta (\c)\vert,$ 
\item $\hat \zeta_{\epsilon}(t)\geq 0$ is decreasing for $\a\leq t\leq\frac{\a+\c}{2}$ and increasing for $\frac{\a+\c}{2}\leq t\leq \b$,
\item $\hat \zeta_{\epsilon}(t)$ converges uniformly to $0$ for $\epsilon\rightarrow 0^+$ on every compact subset of $[\a, \b),$
\item $\hat \zeta_{\epsilon}(t)=\vert \eta (\b)-\eta (\c)\vert{O}\left(e^{\sqrt{\frac{m}{\epsilon}}\hat\chi (t)}\right)$ where $\hat\chi (t)=t-\b$ for $\frac{\a+\c}{2}\leq t\leq \b$ and  $\hat\chi (t)=\c-\b+\a-t$ for $\a\leq t<\frac{\a+\c}{2}.$
\end{enumerate}

\centerline{} 

The correction function 
$$\vc=-\frac{\left(\psi_\epsilon (\a)-\psi_\epsilon(\c)\right)}{(\eta (\c)-\eta (\a))}\zeta_\epsilon(t)+\frac{\left(\psi_\epsilon (\c)-\psi_\epsilon(\b)\right)}{\vert \eta (\b)-\eta (\c)\vert}\hat \zeta_\epsilon(t)+\psi_\epsilon(t)$$
 where
{\setlength\arraycolsep{1pt}
\begin{eqnarray*} 
 \psi_\epsilon(t)&=&\frac{\lambda\vert \eta (\c)-\eta (\a)\vert }{D\sqrt{m\epsilon}}t\Big(e^{\sqrt{\frac m{\epsilon}}(\b-t)}+e^{\sqrt{\frac m{\epsilon}}(t-\b)} \\
&-&e^{\sqrt{\frac m{\epsilon}}(\c-t)}-e^{\sqrt{\frac m{\epsilon}}(t-\c)}\Big)
\end{eqnarray*}}
converges uniformly to $0^+$ on $[\a,\b]$ for $\epsilon\rightarrow 0^+.$

\centerline{}

 Theorem \ref{maintheorem} implies that $y_\epsilon(t)=\eta (t)+O(\epsilon)$ on every compact subset of $(\a,\b)$ and  $$\lim\limits_{\epsilon\rightarrow0^+ }y_\epsilon(\a)=\lim\limits_{\epsilon\rightarrow0^+ }y_\epsilon(\b)=\lim\limits_{\epsilon\rightarrow0^+ }y_\epsilon(\c)=\eta (\c).$$
Consequently,
$$\lim\limits_{\epsilon\rightarrow0^+ }g\left(y_\epsilon(\a)\right)=\lim\limits_{\epsilon\rightarrow0^+ }g\left(y_\epsilon(\b)\right)=\lim\limits_{\epsilon\rightarrow0^+ }g\left(y_\epsilon(\c)\right)=g\left(\eta (\c)\right).$$
Due to the assumption that $g$ is strictly monotone, the boundary layer effect occurs at the point $\a$ or/and $\b$ in the case when $\eta (\a)\neq \eta (\c)$ or/and  $\eta (\b)\neq \eta (\c).$ 

\section{Approximation of realization of SPS}

The application of numerical methods may give rise to difficulties when the singular perturbation parameter $\epsilon$ tends to zero, especially in the nonlinear case. Then the mesh needs to be refined substantially to grasp the solution within the boundary layers (piecewise uniform mesh of Shishkin-type; see, e.g. \cite{MRS}, \cite{RiQu} and the references therein).  The advantage of our approach is that we have to solve only on the parameter $\epsilon$ independent limiting problem $ky=f\left(u(t),y\right),$ see the assumption (A1).
Then a singular perturbation method is applied to obtain an approximate solution of SPS (\ref{def_DE2}), (\ref{def_BC2}) composed of a solution $\eta$ of reduced problem, small constant and two boundary layer functions to recover the lost nonlocal boundary conditions in the degeneration process.

We use the linear combination of the functions $\eta(t), \zeta_{\epsilon}(t)$ and $\hat\zeta_{\epsilon}(t)$ to approximate the exact solution of SPS (\ref{def_DE2}), (\ref{def_BC2}) by the following way. For $\eta \left(\b\right)-\eta \left(\c\right)\leq0$  we define the approximate realization $\tilde y_\epsilon(t)$ of SPS (\ref{def_DE2}), (\ref{def_BC2}) by
\begin{equation}\label{approx1}
\tilde y_\epsilon(t)=\eta (t)+\zeta_\epsilon(t)+\hat \zeta_\epsilon(t)+C\epsilon 
\end{equation}
and analogously, for $\eta \left(\b\right)-\eta \left(\c\right)\geq0$  we define 
\begin{equation}
\tilde y_\epsilon(t)=\eta (t)+\zeta_\epsilon(t)-\hat \zeta_\epsilon(t)-C\epsilon
\end{equation}
where the $\epsilon-$independent constant $C$ is defined in Theorem \ref{maintheorem}.

It is not difficult to verify that $\tilde y_\epsilon(t)$ satisfies the boundary conditions (\ref{def_BC2}) and
\begin{equation*}
\lim\limits_{\epsilon\rightarrow0^+ }\tilde y_\epsilon(\a)=\eta (\c)=\lim\limits_{\epsilon\rightarrow0^+ }\tilde y_\epsilon(\b).
\end{equation*}
Further,
\begin{enumerate}
\item for $\eta \left(\b\right)-\eta \left(\c\right)\leq0$ and $\eta \left(\c\right)-\eta \left(\a\right)\leq0$ we obtain the inequality 
\begin{equation}\label{inequality1}
-\vc\leq\tilde y_\epsilon(t)-y_\epsilon(t)\leq 2\hat \zeta_\epsilon(t)+2C\epsilon,
\end{equation}
\item for $\eta \left(\b\right)-\eta \left(\c\right)\geq0$ and $\eta \left(\c\right)-\eta \left(\a\right)\geq0$ 
\begin{equation}\label{inequality2}
-\vc\leq y_\epsilon(t)-\tilde y_\epsilon(t)\leq 2\hat \zeta_\epsilon(t)+2C\epsilon,
\end{equation}
\item for $\eta \left(\b\right)-\eta \left(\c\right)\leq0$ and $\eta \left(\c\right)-\eta \left(\a\right)\geq0$ 
\begin{equation}\label{inequality3}
0\leq\tilde y_\epsilon(t)-y_\epsilon(t)\leq\vc+2\hat \zeta_\epsilon(t)+2C\epsilon,
\end{equation}
\item for $\eta \left(\b\right)-\eta \left(\c\right)\geq0$ and $\eta \left(\c\right)-\eta \left(\a\right)\leq0$ 
\begin{equation}\label{inequality4}
0\leq y_\epsilon(t)-\tilde y_\epsilon(t)\leq\vc+2\hat \zeta_\epsilon(t)+2C\epsilon.
\end{equation}
\end{enumerate}

The right sides of the inequalities (\ref{inequality1})--(\ref{inequality4}) are $O(\epsilon)$ on every compact subset of $[ \a, \b).$
On the other hand, taking into consideration the facts that $\tilde y_\epsilon(\a)=\tilde y_\epsilon(\b),$  $y_\epsilon(\a)=y_\epsilon(\b)$ and monotonicity of the functions $\vc+2\hat \zeta_\epsilon(t)+2C\epsilon$ and $2\hat \zeta_\epsilon(t)+2C\epsilon$ with respect to the variable $t$ in a left neighbourhood of $\b$ for small $\epsilon,$ we have
\begin{equation*}
\left\vert y_\epsilon(t)-\tilde y_\epsilon(t)\right\vert\leq O(\epsilon)
\end{equation*}
on $[ \a, \b],$ that is, $\tilde y_\epsilon(t)$ is $O(\epsilon)$ accurate approximation of exact solution $y_\epsilon(t)$ of (\ref{def_DE2}), (\ref{def_BC2}) on the whole interval $[ \a, \b].$ We also see that $\left\vert\tilde w_\epsilon(\a)\right\vert\rightarrow \infty$ and $\left\vert\tilde w_\epsilon(\b)\right\vert\rightarrow \infty$ for $\epsilon\rightarrow 0^+,$ where $\tilde w_\epsilon\equiv \tilde y'_\epsilon.$ Thus, $\tilde y_\epsilon(t)$ is a good approximation of the boundary layers arising in the endpoints of the considered interval $[ \a, \b].$

We remark that in the special case when $C=0,$ that is, if $\eta$ is a first-degree polynomial function or a piecewise linear function (in the second case a small generalization of Theorem \ref{maintheorem} is needed) we obtain the exponential convergence rate of $\tilde y_\epsilon$ to $y_\epsilon$ on $[ \a, \b]$ for $\epsilon\rightarrow 0^+.$

We remind, that $\tilde y_\epsilon(t)=\eta(t)$ is not an appropriate approximation of $ y_\epsilon(t)$ because do not respect the possible appearance of boundary layers.

\centerline{}
 
Consider SPS with quadratic nonlinearity of the form
\begin{equation}\label{example}
\epsilon y''+ky=y^2+u(t),\quad k<0,\quad u\in C^2\left([\a,\b]\right) 
\end{equation}
with the boundary conditions (\ref{def_BC}). The assumptions of Theorem \ref{maintheorem} are satisfied  if and only if there exists $\lambda>0$ such that
{\setlength\arraycolsep{1pt}
\begin{eqnarray}
\frac14\left(k^2-(\lambda-k)^2\right)&<& u(t)<\frac14\left(k^2-(\lambda+k)^2\right)\quad\mathrm{on}\quad[ \a,\b]\label{c1}\\
\left\vert u(\c)-u(\a)\right\vert&<&\frac18\left(\lambda-k-\iota (\a)\right)\left(\iota(\a)+\iota(\c)\right)\label{c2}\\
\left\vert u(\b)-u(\c)\right\vert&<&\frac18\left(\lambda-k-\iota(\b)\right)\left(\iota(\b)+\iota(\c)\right)\label{c3}\\
\left\vert u(\c)-u(\a)\right\vert&<&\frac18\left(\lambda+k+\iota(\a)\right)\left(\iota(\a)+\iota(\c)\right)\label{c4}\\
\left\vert u(\b)-u(\c)\right\vert&<&\frac18\left(\lambda+k+\iota(\b)\right)\left(\iota(\b)+\iota(\c)\right),\label{c5} 
\end{eqnarray}}
where $\iota(t)=\st.$ 

\begin{figure}[htb]
\begin{center}
\includegraphics[width=0.75\textwidth, keepaspectratio]{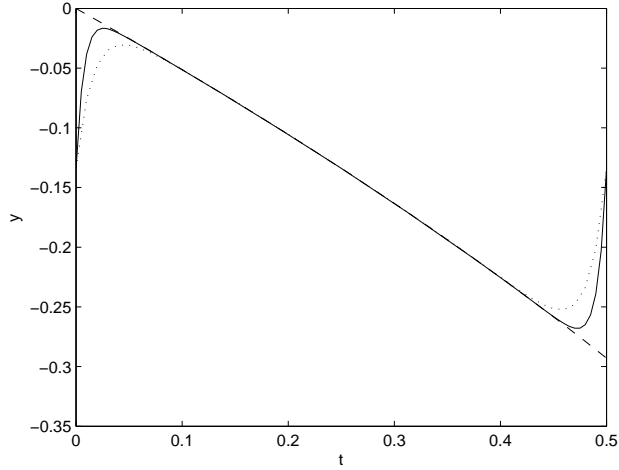}
\caption{Boundary layer phenomenon for solution of singularly perturbed problem $\epsilon y''-2y=y^2+t,$  $y(0)=y(1/4)=y(1/2)$ (the solid line) with  $\epsilon=0.0001.$    The dotted and dashed lines represent the approximate solution $\tilde y_\epsilon(t)$ (with $\lambda=1.6$) and solution of reduced problem, the function $\eta(t)=-1+\sqrt{1-t},$ respectively}
\label{solution}
\label{Figure1}
\end{center}
\end{figure}

For an illustrative example let we consider the problem (\ref{example}), (\ref{def_BC}) with $k=-2,$ $u(t)=t,$ $\a=0,$ $\b=1/2,$ $\c=1/4$ and $g=\mathrm{id}.$ It is not difficult to verify that the solution $\eta(t)=-1+\sqrt{1-t}$ of reduced problem satisfies the conditions (\ref{c1})--(\ref{c5})  for every $\lambda\in\left(\frac2{\sqrt2+\sqrt3}+2-\sqrt2,2\right).$ Thus, on the basis of Theorem \ref{maintheorem}, there exists $\epsilon_0=\epsilon_0(\lambda)$ such that for every  $\epsilon \in (0,\epsilon_0]$ the problem $\epsilon y''-2y=y^2+t,$  (\ref{def_BC}) has in ${H}(\eta)$ the unique solution which is $O(\epsilon)$ close to the approximate solution (\ref{approx1}) on $[ \a, \b]$ (Fig. \ref{Figure1}), that is, to the function
\begin{equation*}
\tilde y_\epsilon(t)=-1+\sqrt{1-t}+\zeta_\epsilon(t)+\hat \zeta_\epsilon(t)+\epsilon\left[(2-\lambda)\sqrt2\right]^{-1}.
\end{equation*}

%\iffalse
%\fi
In the context of previous analysis of the steady--state solutions of 1-D heat transfer equation, it would be interesting to investigate the occurrence of boundary layers for $\epsilon\rightarrow 0^+$ of perturbed, non-stationary 1-D heat transfer equation, written in the usual form as 
\begin{equation*}
\frac{\partial y}{\partial t}=\epsilon\frac{\partial^2 y}{\partial x^2}+ky-f\left(u(x),y\right)  
\end{equation*}
subject to the nonlocal boundary conditions 
\[
v(x_i,t)=v(x_m,t)=v(x_f,t),\quad x_i<x_m<x_f,\quad t\in [0,\infty),
\]
where $v(x,t)=g(y(x,t)).$ The solution $y_\epsilon (x,t)$ represents the temperature at point $x$ of the heated bar in the time $t,$ $x\in [ x_i,x_f],$ $t\in[ 0,\infty).$ For the initial value problems, the numerical analysis of non-stationary reaction-diffusion systems shows on the presence of boundary layer phenomenon (see, e.g. \cite{Shi}).

\section{Feedback control of semilinear SPS}

In this section we consider SPS (\ref{def_system1}), (\ref{def_system22}), (\ref{def_system3}) with
\begin{equation}\label{def_system22}
\epsilon w'(t)=-ky(t)+f(y(t))+u(t).
\end{equation}
Let
\begin{equation*} 
\left\vert f'(y)\right\vert\leq\lambda<-k
\end{equation*}
for $y\in \mathbb{R}.$ Moreover, assume that $g\in C^1$ and $g_{-1}\in C^2$ on $\mathbb{R}$ where $g_{-1}$ denotes an inverse function for $g.$

Now, if $v^0\in C^2\left([\a,\b]\right)$ is desired output of SPS (\ref{def_system1}), (\ref{def_system22}), (\ref{def_system3}) satisfying (\ref{def_BC}) then it is easy to verify that an adequate feedback control input $u^0$ to obtain close $v^0$ output is
\begin{equation*}\label{input}
u^0(t)=kg_{-1}\left(v^0(t)\right)-f\left(g_{-1}\left(v^0(t)\right)\right).
\end{equation*}
Hence $\eta^0(t)=g_{-1}\left(v^0(t)\right)$ and an observable realization  $g\left(y^0_\epsilon\right)$ of system (\ref{def_system1}), (\ref{def_system22}), (\ref{def_system3}) with the boundary condition (\ref{def_BC})  is $O(\epsilon)$ close to the 
$g\left(\tilde y^0_\epsilon(t)\right).$ 

Indeed, as follows from the Lagrange Theorem and (\ref{inequality1})--(\ref{inequality4}),
{\setlength\arraycolsep{1pt}
\begin{eqnarray*}
\left\vert g\left(y^0_\epsilon(t)\right)-g\left(\tilde y^0_\epsilon(t)\right)\right\vert&\leq& \mu\left\vert y^0_\epsilon(t)-\eta^0(t) \right\vert\\
&\leq& \mu\frac{\epsilon}{m}\max\left\{\left\vert \eta^{0''}(t)\right\vert; t\in[ \a,\b]\right\}
\end{eqnarray*}}
where $\mu=\max\left\{\left\vert g'(y)\right\vert; (t,y)\in H\left(\eta^0\right)\right\}.$

\section{Unsolved controllability problem }

Consider the dynamical model described by singularly perturbed differential equation 
\begin{equation}\label{def_DE_unsolved}
\epsilon y''(t)+\frac 12\tilde f\left(u(t),y(t)\right)=0,
\end{equation} 
where $\tilde f=2(ky-f)\in C\left(\mathbb{R}^{2}\right)$ (see (\ref{def_DE2})), $u\in C\left([0,\b]\right)$ is a continuous control input and $0<\epsilon<<1$ is a singular perturbation parameter. Let $\tilde f\neq 0,$ and without loss of generality we will assume that $\tilde f>0$ and $\a=0.$ In this case, the reduced problem $\tilde f\left(u(t),y(t)\right)=0$ does not have a solution $\eta$ (Assumption (A1)), which was the crucial assumption to prove Theorem \ref{maintheorem}. 

Denote by $\{t^*_{i,\epsilon}\}$ the set of turning points in $(0,\c)$ of exact solutions $y_\epsilon$ for problem (\ref{def_DE_unsolved}) satisfying $y_\epsilon(0)=y_\epsilon(\c),$ that is, $y_\epsilon'(t^*_{i,\epsilon})=0$ and $y_\epsilon''(t^*_{i,\epsilon})\neq0.$ For the problems considered in the previous sections, the turning points are determined for small $\epsilon$ with sufficient precision by the turning points of the solution $\eta$ of reduced problem. Obviously, for (\ref{def_DE_unsolved}) there is only one turning point $t^*_\epsilon$ of the solution $y_\epsilon$ on $[0,\b],$ and in $t^*_\epsilon$ acquires its local and global maximum on $[0,\c]$ and it is possible to steer the control system (\ref{def_DE_unsolved}) from the state $y_\epsilon(0)$ to the state $y_\epsilon(\c),$ $0<\c<\b,$ satisfying $y_\epsilon(0)=y_\epsilon(\c)$ with an arbitrary second boundary condition and for every small $\epsilon.$

Now we will analyze the location of this turning point.
 
Let consider a special case of (\ref{def_DE_unsolved}) when  $\tilde f\left(u(t),y(t)\right)\equiv \tilde f\left(u_0,y(t)\right),$ that is, the nonlinear mathematical model 
\begin{equation}\label{def_DE_auton}
\epsilon y''(t)+\frac12 \tilde f\left(u_0,y(t)\right)=0,
\end{equation} 
with the initial conditions $y_\epsilon(0)=y_{0,\epsilon},$ $y'_\epsilon(0)=y_{1,\epsilon},$ where $y_{0,\epsilon},y_{1,\epsilon}$ are the arbitrary real numbers. Obviously, $y_{1,\epsilon}>0,$ because in the case $y_{1,\epsilon}\leq 0$ the solution $y_\epsilon$ of (\ref{def_DE_auton}) satisfying $y_\epsilon(0)=y_\epsilon(\c)$ has a local minimum at some $t_0\in(0,\c)$ with $y''_\epsilon(t_0)\geq 0$ which contradicts to the assumption on positivity of the function $\tilde f.$ Denote by $\tilde F_{u_0}$ the antiderivative of $\tilde f(u_0,y),$ that is, $\tilde F_{u_0}=\int \tilde f(u_0,y)\dy.$ The function $\tilde F_{u_0}$ is strictly increasing and by $\tilde F^{-1}_{u_0}$ we denote an inverse function to $\tilde F_{u_0}.$ Integrating the differential equation (\ref{def_DE_auton}) we have
\begin{equation}\label{integrating}
\epsilon (y'_\epsilon(t))^2+\tilde F_{u_0}(y(t))=\epsilon y_{1,\epsilon}^2+\tilde F_{u_0}(y_{0,\epsilon}).
\end{equation}

Now applying the standard methods we obtain that for every $t\in[0,\b],$ $y_\epsilon(t)$ is an unique root of the equation
\begin{equation}\label{solution_auton}
\pm2\epsilon\int\limits_{\sqrt{\epsilon y_{1,\epsilon}^2-\int\limits_{y_{0,\epsilon}}^{y_\epsilon(t)}\tilde f(u_0,s)\ds }}^{\sqrt{\epsilon y_{1,\epsilon}^2}}\left[\tilde f\left(\tilde F^{-1}_{u_0}\left(\tilde F_{u_0}(y_{0,\epsilon})+\epsilon y_{1,\epsilon}^2-z^2\right)\right) \right]^{-1}\dz=t,
\end{equation} 
where the sign $+(-)$ on the subintervals of $[0,\b]$ with $y'_\epsilon\geq 0$ ($y'_\epsilon<0$), that is, for $t\in(0,t^*_\epsilon]$ ($t\in(t^*_\epsilon,\b]$) is considered, respectively. 

Taking into consideration that $y_\epsilon'(t^*_\epsilon)=0$ we have 
\begin{equation}\label{in_turning_point}
\tilde F_{u_0}(y(t^*_\epsilon))=\epsilon y_{1,\epsilon}^2+\tilde F_{u_0}(y_{0,\epsilon}).
\end{equation}
Thus for computation of the turning point we obtain from (\ref{solution_auton}) the equation 
\begin{equation*}%\label{turning_point}
2\epsilon\int\limits_{0}^{\sqrt{\epsilon y_{1,\epsilon}^2}}\left[\tilde f\left(\tilde F^{-1}_{u_0}\left(\tilde F_{u_0}(y_{0,\epsilon})+\epsilon y_{1,\epsilon}^2-z^2\right)\right) \right]^{-1}\dz=t^*_\epsilon.
\end{equation*} 

To illustrate this theory, let us consider (\ref{def_DE_auton})  with $\tilde f\left(u_0,y(t)\right)=e^y.$ The solution of initial problem is
\begin{equation}\label{exact_solution}
y_\epsilon (t)=\ln\left[c_1-c_1\left(\frac{e^{\mp\frac{\sqrt{c_1}}{\epsilon}(t+c_2)}-1}{e^{\mp\frac{\sqrt{c_1}}{\epsilon}(t+c_2)}+1}\right)^2\right],
\end{equation}
where the sign $-(+)$ on the subintervals of $[0,\b]$ with $y'_\epsilon\geq0$ ($y'_\epsilon<0$) holds, respectively. The constants $c_1,$ $c_2$ are
\begin{equation*}
c_1=\epsilon y_{1,\epsilon}^2+\tilde F_{u_0}(y_{0,\epsilon}),\quad c_2=-\frac{\epsilon}{\sqrt{c_1}}\ln\frac{\sqrt{c_1}+\sqrt\epsilon y_{1,\epsilon}}{\sqrt{c_1}-\sqrt\epsilon y_{1,\epsilon}}.
\end{equation*}
From (\ref{in_turning_point}) we have $y_\epsilon(t^*_\epsilon)=\ln c_1.$ Thus, as follows from  (\ref{exact_solution}), $t^*_\epsilon+c_2=0$ and we obtain
\begin{equation}\label{turning_point_ev}
t^*_\epsilon=\frac{\epsilon}{\sqrt{c_1}}\ln\frac{\sqrt{c_1}+\sqrt\epsilon y_{1,\epsilon}}{\sqrt{c_1}-\sqrt\epsilon y_{1,\epsilon}}.
\end{equation}

On the other hand, from (\ref{exact_solution}), equating $y_\epsilon(0)$ and $y_\epsilon(\c)$  we get $2c_2+\c=0.$ Comparing this with (\ref{turning_point_ev}) we obtain
\begin{equation*}
t^*_\epsilon=\frac \c2.
\end{equation*}

The following questions arise in this context:  
\begin{itemize}
\item[(i)]Where is located the turning point $t^*_\epsilon$ for nonlinear singularly perturbed system (\ref{def_DE_unsolved}) with $\tilde f>0$ subject to required boundary condition $y_\epsilon(0)=y_\epsilon(\c),$ $0<\c<\b$ in general? Does have the position independent of singular perturbation parameter $\epsilon$?  
\item[(ii)]Can be controlled a location of turning point by using an appropriate control signal $u$?   
\end{itemize}

\section*{Acknowledgments}
I would like to express our gratitude to the referees for all the valuable and constructive comments.

\end{document}